\documentclass[11pt,a4paper,twoside]{amsart}
\usepackage{amsmath}
\usepackage{amssymb}
\usepackage{latexsym}

\newcommand{\co}{\colon\thinspace}    
\newcommand{\fnote}[1]{\footnote{\small sharp1}}

\newcommand{\N}{{\mathbb N}}
\newcommand{\Z}{{\mathbb Z}}
\newcommand{\R}{{\mathbb R}}
\newcommand{\Q}{{\mathbb Q}}

\newcommand{\T}{{\mathbb T}}
\newcommand{\Sbar}{{\mathbb S}}

\newcommand{\rat}{{\rm rat}}

\newtheorem{theorem}{Theorem}[section]
\newtheorem{proposition}[theorem]{Proposition}

\newtheorem{definition}[theorem]{Definition}
\newtheorem{lemma}[theorem]{Lemma}

\title{ Ma\~{n}\'e's conjectures in codimension one}
\author{Ugo Bessi, Daniel Massart}
\date{\today}
\begin{document}

\begin{abstract}
We prove Ma\~n\'e's conjectures (\cite{Mane96}) in the context of codimension one Aubry-Mather theory.
\end{abstract}
\maketitle

\section{Introduction}
We study variational problems on tori in the spirit of \cite{Moser}. The objects we are interested in are maps $u$ from $\R^n$ to $\R$ which minimize globally the integral 
	\begin{equation}\label{integral_to_be_minimized}
	\int_{\R^n}L(x,u,\nabla(u) )dx
\end{equation}
where the cost function $L$ is called the Lagrangian of the problem. 
This theory is also known as codimension one Aubry-Mather theory, because it generalizes the classical Aubry-Mather theory of twist maps. It runs parallel to the ''dimension one'' theory of Bangert \cite{Bangert90}, Mather \cite{Mather91}, Ma\~n\'e \cite{Mane96} and Fathi \cite{Fathi}. 

We begin by recalling the hypotheses on the Lagrangian. Let \newline
$L(x_1,\dots,x_n,u,p_1,\dots,p_n)$ be a Lagrangian such that

\begin{itemize}
\item  (H1) : $L\in C^{l,\gamma}(\R^{2n+1})$, $l\ge 2$, $\gamma> 0$.

\item (H2) : $L$ has period $1$ in $x_1,\dots,x_n,u$.

\item (H3) : There is $\delta>0$ such that
$$\delta I\le\frac{\partial^2L}{\partial p_i\partial 
p_j}\le{\frac{1}{\delta}}I$$
where $I$ denotes the identity matrix on $\R^n$.

\item (H4) : There is $C>0$ such that
$$\left\vert\frac{\partial^2L}{\partial p\partial x}\right\vert+
\left\vert\frac{\partial^2L}{\partial p\partial u}\right\vert
\le C(1+|p|)$$
$$\left\vert\frac{\partial^2L}{\partial x\partial x}\right\vert+
\left\vert\frac{\partial^2L}{\partial u\partial x}\right\vert+
\left\vert\frac{\partial^2L}{\partial u\partial u}\right\vert
\le C(1+|p|^2) .  $$

\end{itemize}
The main example we have in mind is a Lagrangian of the form
	\[ L(x,u, \nabla u)= \frac{1}{2}\left|\nabla u(x)\right|^2 +f(x,u)
\]
where $f \in C^{l,\gamma}(\R^{n+1})$ is $\Z^{n+1}$-periodic. Observe that for any Lagrangian $L$ satisfying Hypothesis (H1-4) and for any $\Z^{n+1}$-periodic 
$f \in C^{l,\gamma}(\R^{n+1})$, $L+f$ is again a Lagrangian  satisfying Hypothesis (H1-4). Adding a function to the Lagrangian is also called perturbing the Lagrangian by a potential. In this paper, after Ma\~n\'e (\cite{Mane95}),  the phrase "for a generic Lagrangian $L$, Property P holds" means "for any   Lagrangian $L$, there exists a residual subset $\mathcal{O}(L)$ of the set of potentials, such that for any $f$ in $\mathcal{O}(L)$, Property P holds for $L+f$".

Since the integral (\ref{integral_to_be_minimized}) is infinite in general, we must explain what we mean by minimizing in (\ref{integral_to_be_minimized}).
We say that $u\in W^{1,2}_{loc}(\R^n)$ is a minimizer for $L$ if
\begin{equation}\label{compactly_supported_variation}
\int_{\R^n}[L(x,u+\phi,\nabla(u+\phi) )-L(x,u,\nabla u)]dx\ge 0
\qquad \forall \phi\in C^\infty_0(\R^n).
\end{equation}
Since $L$ is periodic, if $u$ is a minimizer and 
$(k,j)\in\Z^n\times\Z$, then $u(x+k)+j$ is a minimizer,  too; we say 
that $u$ is non self intersecting if 
\begin{equation}\label{non_self-intersecting}
\forall(k,j)\in\Z^n\times\Z,\quad\hbox{either}\quad u(x+k)+j>u(x)\ \forall x$$
$$\quad\hbox{or}\quad u(x+k)+j<u(x)\ \forall x\quad\quad\hbox{or}\quad 
u(x+k)+j=u(x)\ \forall x  .
\end{equation}
In \cite{Moser}, it is proven that non self-intersecting minimizers  lie within finite distance of some hyperplane:

\begin{theorem}[\cite{Moser}]\label{uniform bounds on minimizers} Let 
$u\in W^{1,2}_{loc}(\R^n)$ be minimal and non self intersecting; then there exists $\rho\in\R^n$ and a constant $C_L(||\rho||)>0$, depending only on $L$ and $||\rho||$, such that
$$||u-u(0)-\rho\cdot x||_{C^{l,\gamma}(\R^n)}\le C_L(||\rho||) . $$

\end{theorem}
In particular, any minimizer $u\in W^{1,2}_{loc}(\R^n)$ is actually as regular as the Lagrangian.
The vector $\rho$ is called the rotation vector, or the slope, of $u$; an important fact is that there are minimal, non self intersecting solutions of any rotation vector. 
\begin{theorem}[\cite{Moser}] 
For any $\rho\in\R^n$, there is a minimal, non self intersecting solution of slope $\rho$.
\end{theorem}

\begin{definition}
A minimal, non self intersecting solution of slope $\rho$, is called a $(L,\rho)$-minimizer. When $\rho\in\Q^n$, we can consider the subclass of periodic minimizers: we say that a $(L,\rho)$ minimizer $u$ is periodic if $u(x+k)+j=u(x)$ for all $(k,j)\in\Z^n\times\Z$ such that 
$\rho\cdot k+j=0$. If $u$ is a  $(L,\rho)$-minimizer, with  
$\rho \in \Q^n$, then $u$ is either periodic, or asymptotic to some periodic $(L,\rho)$-minimizer (see \cite{Bangert}).
\end{definition}

We want to study uniqueness of $(L,\rho)$-minimizers; since we saw before that, if $u$ is a $(L,\rho)$-minimizer, also 
$u(\cdot+k)+j$ is such, we have to identify $u$ with its integer translations. Even with this identification, the answer is negative, because in \cite{Bangert} it is proven that, if $\rho\not\in\Q^n$, or if $\rho\in\Q^n$ and $n\ge 2$, there are always uncountably many 
$(L,\rho)$-minimizers. The situation changes if we look at the currents induced by minimizers (see section \ref{Action of a current} for the precise definitions). Indeed, we are able to prove that, generically, all $(L,\rho)$-minimizers induce the same current; with the added bonus that, if $\rho$ is irrational, we can drop the "generically."

The problem of uniqueness can be formulated not only for 
$(L,\rho)$-minimizers, but also for the dual notion of $(L-c)$-minimizers. We briefly explain what we mean; we recall, that, as proven in \cite{Senn91}, a mean action is defined.

\begin{theorem}[\cite{Senn91}]\label{Senn1}
For any $(L,\rho)$-minimizer $u$, the following limit exists and depends only on $L$ and $\rho$: 
$$\lim_{R\to\infty}{\frac{1}{|B(0,R)|}}
\int_{B(0,R)}L(x,u,\nabla u)dx :=\beta(\rho) . $$
Moreover, the function $\beta$ is strictly convex and superlinear.

\end{theorem}

Since $\beta$ is strictly convex, its Legendre-Fenchel transform, traditionally denoted by $\alpha$, is $C^1$; it is easy to see that 
$-\alpha(c)$ is the minimum, over all $u$ minimal and non self intersecting, of 
$$\lim_{R\to\infty}{\frac{1}{|B(0,R)|}}
\int_{B(0,R)}[L(x,u,\nabla u)-c\cdot \nabla u]dx   . $$
Note that for any $c$ in $\R^n$, the Lagrangian $L(x,u,\nabla u)-c\cdot \nabla u$, denoted $L-c$ for short, still satisfies Hypothesis (H1-4). A minimal, non self-intersecting $u$ such that  
\[
\lim_{R\to\infty}{\frac{1}{|B(0,R)|}}
\int_{B(0,R)}[L(x,u,\nabla u)-c\cdot \nabla u]dx   =  -\alpha(c)
\]
is called a $(L-c)$-minimizer. As for $(L,\rho)$-minimizers, we may ask about the uniqueness of the $(L-c)$-minimizer for a given $c$, and similarly the question should be rephrased in terms of currents. One difference between $(L,\rho)$-minimizers and $(L-c)$-minimizers is that we don't know a priori when an $(L-c)$-minimizer is periodic, so another question we adress is \textit{how large is the set of $c$ for which $(L-c)$-minimizers have a rational slope ?} Note that by Fenchel duality an $(L-c)$-minimizer is an $(L,\alpha'(c))$-minimizer, so the question boils down to \textit{how large is the set of $c$ for which $\alpha'(c) \in \Q^n$ ?}.

Now we can state our result.
\begin{theorem}
For a generic Lagrangian satisfying Hypothesis (H1-4),
\begin{itemize}
	\item for every $\rho \in \R^n$, the $(L,\rho)$-minimizers induce a unique current; if $\rho$ is rational, there is a unique periodic 
	$(L,\rho)$-minimizer
	\item for every $c \in \R^n$, the $(L-c)$-minimizers induce a unique current
	\item there exists an open dense subset $U$ of $\R^n$ such that for every $c \in U$, $\alpha'(c) \in \Q^n$.
\end{itemize}
\end{theorem}

Our theorem solves, in the affirmative, the codimension one versions of the problems posed by Ma\~n\'e in \cite{Mane95}, \cite{Mane96}. In the ''dimension one'' theory much less is known. The best result about the first point of the theorem is that of \cite{Bernard-Contreras}, which says that for a generic Lagrangian $L$ on a manifold of dimension $n$, for every homology class $\rho$, there exists at most $n+1$ 
$(L,\rho)$-minimizing currents. 
The second point of the theorem  is trivially false in the dimension one theoretical setting (see Hedlund's example in \cite{Bangert90}). To be precise about the third point, recall that the problem originally proposed by Ma\~n\'e was : \textit{is it true that for a generic Lagrangian $L$, there exists a dense open subset $U$ of the cohomology of the configuration space such that for any $c \in U$, there exists a unique minimizing measure, and it is supported on a periodic orbit}. This   is true, by \cite{Osvaldo}, when the base manifold is the circle and the Lagrangian depends periodically on time, and by \cite{ijm} when the base manifold has dimension two and the Lagrangian is autonomous. In the codimension one theory,  the notion corresponding to \textit{minimizing measure} is that of \textit{recurrent minimizer}. Thus, in this case Ma\~n\'e conjecture follows by the first and third points of theorem 1.5.

Thus Ma\~n\'e's conjectures seem taylor-made for the codimension one case. One possible reason for this is that Ma\~n\'e had in mind the twist map case, which in some respects is more typical of the codimension one case than it is of the dimension one case.

\textbf{Acknowledgements} The second author was partially supported by the ANR project ''Hamilton-Jacobi et th\'eorie KAM faible''.
\section{The derivative of $\alpha$ is rational on a dense set}
We define
$$\rat(\rho,1)=\mbox{Vect} \left((\rho,1)^\perp\cap\Z^{n+1}\right)  $$
where $\mbox{Vect}(A)$ denotes the smallest subspace of 
$\R^{n+1}$ containing the set $A$.

Let $\alpha$ and $\beta$ be as in the introduction; we recall that they are dual convex functions; since $\beta$ is superlinear and strictly convex by theorem \ref{Senn1}, $\alpha$ is $C^1$ and superlinear.

We call flat of slope $\rho$ the set 
$$
D_{\rho}=\{ (c,\alpha(c))\quad\colon\quad \alpha^\prime(c)=\rho \}  . 
 $$
We shall need the following result of Senn \cite{Senn95}; it says that the linear space generated by the flat of slope $\rho$ is contained in 
$\rat(\rho,1)$.  If $A\subset\R^{p}$, let 
$L(A)$ be the linear space generated by the differences $a-b$ with 
$a,b\in A$. Clearly, if $0\in A$, then $L(A)=\mbox{Vect}(A)$.

\begin{theorem}[\cite{Senn95}]\label{Senn2} 
Let $L$ be a Lagrangian on $\R^{2n+1}$ satisfying Hypothesis (H1-4), and let
$$D_\rho=\{ (c,\alpha(c))\quad\colon\quad \alpha^\prime(c)=\rho \}  .  $$
Then
$$L(D_\rho)  = \rat(\rho,1) $$ 
unless the recurrent $(L,\rho)$-minimizers (i. e. the periodic ones when $\rho$ is rational, and the functions $u^\alpha$ defined in lemma 5.1 below when 
$\rho$ is irrational) foliate $\T^{n+1}$, in which case $L(D_\rho) = \{0\}$.

\end{theorem}
\begin{proposition}\label{dense}
Let $L$ be a Lagrangian on $\R^{2n+1}$ satisfying Hypothesis (H1-4). Then the set $\left\{c \in \R^n \co \alpha' (c) \in \Q^n \right\}$ is dense in $\R^n$.
\end{proposition}
\proof

Let $U$ be any open subset of $\R^n$.

\textbf{First case} : there exists $c$ in $U$ such that the flat 
$D_{\alpha'(c)}$ of $\alpha$ containing $(c,\alpha(c))$ is reduced to a point. Then, by the convexity of $\alpha$,  
	\[ \forall d \in \R^n\setminus\{ c \},\  
	\left\langle \alpha'(c)-\alpha'(d),c-d \right\rangle > 0.
\]
Let $B$ be a closed ball centered at $c$ and contained in $U$. By Theorem  \ref{Senn1},   $\alpha'$ is continuous.  Hence, by Lemma \ref{appendix_lm1}, $\alpha'(U)$ contains a neighborhood of 
$\alpha'(c)$; thus there exists $d \in  U$ such that $\alpha'(d) \in \Q^n$.

\textbf{Second case} : any $c \in U$ is contained in a non-trivial face of $\alpha$, that is to say, for any $c \in U$, the face $D_{\alpha'(c)}$ of $\alpha$ is not reduced to a point. Then by Theorem \ref{Senn2}, for any $c \in U$, the vector space $L(D_{\alpha'(c)})$  generated by 
$D_{\alpha'(c)}$ is  
$\mbox{rat}(\alpha'(c),1)$, which is a rational subspace of $\R^{n+1}$: it is generated, practically by definition, by integer vectors. There are only countably many rational subspaces of $\R^{n+1}$, so by Baire's theorem (a countable union of nowhere dense subsets of a complete metric space is nowhere dense) there exists an open subset $U_1$ of $U$, and a rational subspace $N_1$ of $\R^{n+1}$, such that for any $c \in U_1$, 
\[
N_1 = L(D_{\alpha'(c)}) = \mbox{rat}(\alpha'(c),1).
\]
Let $M_1$ be the canonical projection to $\R^n$ of $N_1$. If $M_1$ has dimension $n$, then $\alpha^\prime(c)$ is rational; therefore, we shall suppose that $M_1$ is a proper subspace of $\R^n$.

\textbf{Observation} : First let us observe that, if $c \in U_1$, then  
$(c,\alpha(c))$ lies in the relative interior of $D_{\alpha'(c)}$. Indeed, let us take $c \in U_1$, and a convex neighborhood $V$ of $0$ in $M_1$, such that 
$c+V \subset U_1$.  Let us denote by $\tilde\alpha$ the map $\alpha$ restricted to $c+V$; then $\tilde\alpha$ is  affine and convex.
The convexity  is trivial, to prove that $\tilde\alpha$ is  affine, we recall one fact from convex analysis: $\tilde\alpha$ is affine on $c+V$ if and only if, for any $d \in c+V$,  the flat of 
$\tilde\alpha$ containing $(d,\alpha(d))$ has maximal dimension. In our case, $c+V$ is an open set of $c+M_1$, and maximal dimension means the dimension of $M_1$.
Now, $D_{\tilde\alpha'(d)}$ is simply $D_{\alpha'(d)}$ intersected with $(c+V)\times\R$; our assumptions on $D_{\alpha'(d)}$ and 
$V$ yield that
$L(D_{\tilde\alpha'(d)})= \rat(\alpha^\prime(c),1)$, and 
$\rat(\alpha^\prime(c),1)$ has the same dimension as $M_1$.

This proves that $\tilde\alpha$  is affine on  the set 
$c+V$, and that $c+V$ is open in $c+M_1$; in other words, 
$(c,\alpha(c))$ lies in the relative interior of 
$D_{\alpha^\prime(c)}$.

\vskip 1pc

From this we now deduce that for any $c \in U_1$, the map $\alpha$ restricted to $(c+M_1^{\perp}) \cap U_1$, which we denote $\alpha_c$ for simplicity,  is strictly convex at $c$, that is, $c$  is not contained in any non-trivial face of $\alpha_c$. Indeed, \cite{AvsM}, Lemma A.2 says that, if some  $c \in U_1$ is contained in a non-trivial face of 
$\alpha_c$ and the observation above holds, then $c$ is contained in a face $D$ of $\alpha$ such that $L(D)$ properly contains $N_1$; but this contradicts the fact that the flat at $c$ generates $N_1$.

So for any $c \in U_1$, the map  $\alpha_c$ is strictly convex at $c$. Therefore,  by the same argument as in the first case,  for any $c \in U_1$, there exists 
$d \in (c+M_1^{\perp}) \cap U_1$, such that 
$\alpha_c' (d) \in M_1^{\perp} \cap \Q^n$.  

Observe that $M_1^{\perp} \cap \Q^n\not=\{ 0 \}$. To show this, we note that $M_1^\perp\not=\{ 0 \}$, because we are supposing that 
$M_1$ is proper; moreover, $M_1^{\perp} $, being the orthogonal of a rational subspace of $\R^n$, is itself a rational subspace of $\R^n$. 

Now $\alpha' (d)$ is the sum of $\alpha_c' (d)$ and the derivative at $d$ of the restriction of $\alpha$ to $(c+M_1) \cap U_1$, which is the orthogonal projection of $\alpha' (d)$ to $M_1$. The latter lies in $\Q^n \cap M_1$ by Lemma \ref{B2} in the appendix. Therefore $\alpha_c' (d) \in \Q^n $.
\qed

\section{Currents and recurrent minimizers}\label{current}

We  define the current induced by  a minimal $u$. 
For $p\in\R^n$, we  denote by $\omega(x,u)\cdot (p,1)$ the $n$-form 
$\omega$ applied to the $n$-vector
$$\left(
\begin{matrix}
1,&0,&\dots,&0\cr
0,&1,&\dots,&0\cr
\dots,&\dots,&\dots,&\dots\cr
p_1,&p_2,&\dots,&p_n     
\end{matrix}
    \right)    .   $$
    
Let $\omega$ be a smooth $n$-form on the torus and let $u$ be a 
$(L,\rho)$-minimizer; let $B(0,R)$ be the ball of radius $R$ centered at the origin in 
$\R^n$, and let $\left|B(0,R)\right|$ be its Euclidean $n$-dimensional volume. It can be proven that the following limit exists
\begin{equation}\label{corri}
\lim_{R\to\infty}\frac{1}{|B(0,R)|}\int_{B(0,R)}
\omega(x,u(x))\cdot(\nabla u(x),1)dx  .
\end{equation}
We define $T_u(\omega)$ to be the limit above. It is  proven in 
\cite{Bessi09}
that $T_u$ is a $n$-current of finite mass, and that 
$\partial T_u=0$. This means the following: let us denote by 
$\Omega^{(0)}_n$ the set of continuous n-foms on 
$\T^{n+1}$, equipped with the $\sup$ norm; then $T_u$ is a linear, continuous operator on $\Omega^{(0)}_n$ and 
$T_u({\rm d}\eta)=0$ for every $(n-1)$-form $\eta$ of class $C^1$. In particular, we can restrict $T_u$ to the subspace of closed forms and quotient on the exact forms; what we obtain is a linear operator ${\sl slope}_{T_u}$ from $H^n(\T^{n+1})$,
the n-th real cohomology group of $\T^{n+1}$, to $\R$. Thus, 
${\sl slope}_{T_u}$ belongs to the dual of $H^n(\T^{n+1})$, which identifies with the $n$-th homology group $H_n(\T^{n+1})$. On 
$H^n(\T^{n+1})$ we introduce, as a basis, the equivalence classes of the differential forms 
$${{\rm d}\hat x_i}\colon=(-1)^{n-i+1}{\rm d}x_1\wedge\dots\wedge
{\rm d}x_{i-1}\wedge{\rm d}x_{i+1}\wedge\dots\wedge{\rm d}x_{n+1}
\quad 1\le i\le n  $$
$${{\rm d}\hat x_{n+1}}=
{\rm d}x_1\wedge\dots\wedge
{\rm d}x_{n} . $$
On $H_n(\T^{n+1})$ we introduce the basis $e_i$ dual to 
${\rm d}\hat x_i$. 
It is easy to see that, with this choice of the basis, if $u$ is $(L,\rho)$-minimal, then 
${\sl slope}_{T_u}=(\rho,1)$.


Given a current $T$ of finite mass on $\T^{n+1}$, we can define a signed measure $\mu_T$ on $\T^{n+1}$ by the formula
	\begin{equation}\label{mutti}
	T(f{\rm d}x_1\wedge\dots\wedge{\rm d}x_n)=
	\int_{\T^{n+1}}fd\mu_T 
	\end{equation} 
for any function $f$ continuous on the torus.

We note that, by (\ref{corri}) and (\ref{mutti}), if
$T = T_u$, the measure $\mu_T$ is defined by 
	\begin{equation}\label{misura}
	 \int_{\T^{n+1}}f(x,x_{n+1})d\mu_T = \lim_{R \rightarrow\infty}\frac{1}{\left|B(0,R)\right|}\int_{B(0,R)}f(x,u(x))dx  .
	\end{equation}
From the formula above, it is immediate that $\mu_{T}$ is a probability measure.

The following lemma will be useful along the way. 

\begin{lemma}\label{translation, current}
For every $u$ in $M_{\rho}$, for every $(z, z_{n+1})$ in $\Z^n \times \Z$, denoting $v(x) := u(x+z) +  z_{n+1}$, we have $T_u =T_v$.
\end{lemma}

\proof
Take \begin{itemize}
  \item $u$ in $M_{\rho}$
  \item a smooth n-form $ \omega$ on $\T^{n+1}$
  \item $(z, z_{n+1})$ in $\Z^n \times \Z$.
\end{itemize}
We have 
\begin{eqnarray*}
T_u (\omega )  &=&   \lim_{R \rightarrow \infty}\frac{1}{|B(0,R)|}     \int_{B(0,R)}  \omega(x,u(x))\cdot(\nabla u(x),1) dx \\
&=&  \lim_{R \rightarrow \infty}\frac{1}{|B(0,R)|}     \int_{B(-z,R)} \omega(x,u(x+z)+ z_{n+1})\cdot(\nabla u(x+z),1) dx \\
&=&   \lim_{R \rightarrow \infty}\frac{1}{|B(0,R)|}     \int_{B(0,R)} \omega(x,u(x+z)+ z_{n+1})\cdot(\nabla u(x+z),1) dx \\
&=& T_v ( \omega )  .
\end{eqnarray*}
The second equality comes from the change of variables 
$x \mapsto x+z$ and the fact that $\omega$ is $\Z^{n+1}$-periodic, the third one from the fact that $\omega$ is bounded and 
$$
\lim_{R \rightarrow \infty}\frac{ |B(0,R)\setminus B(-z, R)|}{|B(0,R)|} =0.
$$
\qed

\subsection{Action of a current}\label{Action of a current}
This action has been defined for dimension 1 currents in 
\cite{Bernard-Buffoni}; as shown in \cite{Bessi09}, the same definition applies to codimension 1 currents. We are not going to recall this definition here, we only recall some facts; the first one is that this definition extends the notion of mean action we gave in the introduction.
 
Indeed, the following holds: if $u$ is a $(L,\rho)$ minimizer, then the mean action of $T_u$, say $MA(L,T_u)$, is given by
$$MA(L,T_u)=\lim_{R\to\infty}
{1\over{|B(0,R)|}}\int_{B(0,R)}L(x,u,\nabla u)dx = \beta( \rho) . 
$$

\subsection{Rational rotation numbers}\label{rational-currents}

Let 
$\rho\in\Q^n$ and let $r>0$; we define
\begin{equation}\label{group}
\Gamma\colon = \{ k\in\Z^n\quad\colon\quad k\cdot\rho\in\Z \} . 
\end{equation}
It is easy to see that $\Gamma$ is a subgroup of $\Z^n$; actually, it is the projection of ${\rm rat}(\rho,1)\cap \left(\Z^n\times\Z\right)$ to $\Z^n$.
Since 
$\rho$ is rational, $\Gamma$ contains a basis of $\R^n$; in particular, the action of $\Gamma$ on $\R^n$ admits a bounded, measurable fundamental domain $D$.

We define a set $J_r (\rho)$ which will come handy in the next section. The set $J_r (\rho)$ is the set of all functions 
$u\colon\R^n\to\R$  satisfying the three points below:

\begin{itemize}
  \item  $u \in  C^{l,\gamma}(\R^n) \subset W^{1,2}_{loc}(\R^n)$
  \item the $ C^{l,\gamma}(\R^n)$-norm of the map $x \mapsto u(x)-u(0) - \rho\cdot x$ is smaller than $r$.
  \item $u(x + k) + j = u(x)$ whenever $(k, j) \in  \Z^n \times \Z \cap (\rho, 1)^{\bot}$.
\end{itemize}

Then $u$ induces a current $T_u$ by (\ref{corri}). The mean action of $T_u$ is given, as expected, by 
\begin{equation}\label{mean-action}
MA(L,T_u)={1\over |D|}
\int_{D}L(x,u,\nabla u)dx  .
\end{equation}

We shall need a theorem, due to Moser, which says that, when $\rho$ is rational and $r$ is larger than the constant $C_L(||\rho||)$ of theorem \ref{uniform bounds on minimizers}, there are minimizers in the class $J_r (\rho)$. This is in sharp contrast to the dimension one case, where in general there are no periodic minimizers.

\begin{theorem}[\cite{Moser}]\label{periodic} Let $\rho\in\Q^n$, let
$\Gamma$ be defined as in (\ref{group}) and let $D$ be a fundamental domain. We set
$$W=\{ u\in W^{1,2}_{loc}(\R^n)\quad\colon\quad 
u(x+k)-u(x)-\rho\cdot k\equiv 0\quad\forall k\in\Gamma \} . $$
Then
$$\beta(\rho)=\inf\{ 
{1\over{|D|}}\int_DL(x,u,\nabla u)dx\quad\colon\quad u\in W
 \}  .  $$
Moreover, the $\inf$ is a minimum and the functions $u\in W$ on which the minimum is attained are $(L,\rho)$ minimizers. 
 
Conversely, if $u$ is a $(L,\rho)$ minimizer such that 
$u(x+k)-u(x)-\rho\cdot k
\equiv 0$ for any $k\in\Gamma$, then 
 $$\beta(\rho)=
 {1\over{|D|}}\int_DL(x,u,\nabla u)dx   .
 $$

\end{theorem}

\section{Generic uniqueness of periodic minimizers, and of the minimizing currents with rational slope}

Here we prove (Proposition \ref{generic_rho}) that given a rational rotation number $\rho$,  for a generic Lagrangian $L$,  there is   a unique periodic minimizer with rotation number $\rho$. Then we prove (Lemma \ref{unicurrent}) that for such a Lagrangian, all $(L,\rho)$-minimizers, including the non-periodic ones, induce the same current.
\begin{proposition}\label{generic_rho}
Let \begin{itemize}
\item $L$ be a Lagrangian on $\R^{2n+1}$ satisfying Hypothesis (H1-4).
\item $\rho$ be a vector in $\Q^n$.
\end{itemize}
Then there exists a residual subset $\mathcal{O}(L,\rho)$ of $C^{\infty}(\T^{n+1})$ such that for any $f \in \mathcal{O}(L,\rho)$, there is only one periodic  $(L-f,\rho)$-minimizer. Moreover, all 
$(L-f,\rho)$-minimizers induce the same current.
\end{proposition}
\proof 

We will see how this result follows from \cite{Bernard-Contreras}.
This paper considers the following situation:
$$\begin{matrix}
{}&H_r&{}&{}\cr
{}&\cap &{}&{}\cr
E\times&F&{}&{}\cr
{}&\downarrow\pi&\searrow&{}\cr
E\times&G&{\buildrel \alpha \over \longrightarrow}&\R\cr
{}&\cup&{}&{}\cr
{}&K&{}&{}    
\end{matrix}$$
where $E$, $F$, $G$ are topological vector spaces, $\pi$ is a linear continuous map between $F$ and $G$ and $\alpha$ is a bilinear coupling. The hypotheses are the following ones.
\begin{itemize}
	\item The bilinear pairing $\alpha$ is continuous.
	\item $K$ is a compact and convex set, separated by $E$; the 
	latter  means that, if $\eta$ and $\nu$ are two different points
	of $K$, then there exists $u\in E$ such that 
	$\alpha(u,\eta-\nu)\not=0$.
	\item $E$ is a Frechet space.
	\item $H_r$ is compact, convex, and $\pi(H_r)\subset K$.
\end{itemize}

We define $H_r^\ast$ as the set of all the affine, continuous functions on $H_r$; for $L$ in $H_r^\ast$, we  denote by
$MIN_{H_r}(L)$ the set of minima of $L$ over $H_r$.

Under these hypotheses, theorem 5 of \cite{Bernard-Contreras} holds:

\begin{theorem}[\cite{Bernard-Contreras}]\label{BC} For any finite dimensional affine subspace $B$ of $H^\ast_r$, there exists a residual subset 
${\mathcal O}(B)\subset E$ such that, for all $f\in{\mathcal O}(B)$ and 
$L\in B$, we have that $\pi(MIN_{H_r}(L-f))$ is contained in an affine subspace of $G$, whose dimension is not larger than the dimension of $B$.
\end{theorem}

We want to apply this theorem to our situation. To do this, we let
\begin{itemize}
	\item $E$ be the Fr\'echet space $C^{l,\gamma}(\T^{n+1})$ 
	\item $F$ be the space of closed $n$-currents of finite mass on $\T^{n+1}$. This space is the dual of the space $\Omega^{(0)}_n$ of continuous $n$-forms on $\T^{n+1}$, equipped with the $\sup$ norm.
	\item $G$ be the dual space of $C^{0}(\T^{n+1})$, i.e. the space of Borel signed measures on $\T^{n+1}$
	\item $\pi \co F \longrightarrow G$ be the continuous linear map $T \longrightarrow \mu_T$ defined  as in (\ref{mutti}).
	
	\item $\alpha$ be the continuous bilinear pairing between $E$ and $G$ defined by integration
	\item $ K \subset G$ be the metrizable, compact, convex set of Borel probability
measures on $\T^{n+1}$. Observe that $K$ is separated by $E$. 
\end{itemize}

The definition of $H_r$ is a bit trickier. We let $J_r (\rho)$ be the set of functions defined in section 
\ref{rational-currents}.

Define $\tilde{H_r}$ to be the set of
currents of the form $T_u$, with  $u $ in  $J_r (\rho)$, and let $H_r$ be
the weak${}^\ast$ closure of the convex hull of $\tilde{H_r}$. Then  $H_r$ is contained in a ball in $F$, so by the  Banach-Alaoglu Theorem it is compact with respect to the weak$^*$ topology. It is also metrizable because  the space $\Omega_{n}^{0}(\T^{n+1})$ of continuous $n$-forms on 
$\T^{n+1}$, equipped with the $\sup$ norm, is separable. 

We saw in formula (\ref{misura}) that $\pi$ brings any 
$T_u\in\tilde H_r$ to a probability measure, i. e. to an element of $K$; taking convex combinations, the same is true for $H_r$.

We show in lemma \ref{continuity} below that 
$MA(L,\cdot)\in H_r^\ast$, i. e. it is an affine, continuous functional on $H_r$.

Now we can apply Theorem \ref{BC}. For us, $B$ will be a singleton, i. e. 
$B= \{ MA(L,\cdot) \}$.
By theorem \ref{BC} there exists a residual subset 
$\mathcal{O}_r(L)$ of $E$ such that for any 
$f \in \mathcal{O}_r(L)$, $\pi(MIN_{H_r}(L-f))$ is reduced to a point. Clearly, if $r$ is smaller than the constant 
$C_{L-f}(||\rho||)$ of theorem \ref{uniform bounds on minimizers}, the minima in $H_r$ may not correspond to any 
$(L-f,\rho)$-minimizer. That's why we consider 
$$\mathcal{O}(L)=\bigcap_{r\in\N}\mathcal{O}_r(L) . $$
We get that $\mathcal{O}(L)$ is a residual set too and, if 
$f\in\mathcal{O}(L)$, then $\pi(MIN_{H_r}(L-f))$ is reduced to a point for any $r$.  

We show how this implies the thesis.
Let us suppose by contradiction that there are two different 
periodic $(L,\rho)$-minimizers, say $u$ and $v$. We prove below that, if $r> C_{L-f}(||\rho||)$, then $T_u$ and $T_v$ are minimal in 
$H_r$. We recall from \cite{Moser} that the graphs of the two periodic minimals $u$ and $v$ are disjoint; this implies by 
(\ref{misura}) that $\mu_{T_u}$ and $\mu_{T_v}$ are different. In other words, $\pi(MIN_{H_r}(L-f))$ contains at least two elements, while we have just proven that it is reduced to a point; this contradiction proves the thesis. 

Now we show that, if $u$ is a $(L-f,\rho)$-minimizer, then $T_u$ minimizes
$MA(L-f,\cdot)$ in $H_r$ for $r$ large enough. To show this, it suffices to show that the minimum of 
$MA(L-f,\cdot)$ on the currents of $H_r$ coincides with the minimum of $MA(L-f,\cdot)$ on the currents $T_u$, where $u$ is a periodic minimizer. By (\ref{mean-action}) and theorem \ref{periodic}, the latter minimum is $\beta_{L-f}(\rho)$, where by $\beta_{L-f}$ we denote the $\beta$-function of the Lagrangian $L-f$. Thus, it suffices to show that the minimum of $MA(L-f,\cdot)$ on the currents of $H_r$ coincides with 
$\beta_{L-f}(\rho)$. That's what we do next.

We begin by noting the following: let $r>C_{L-f}(||\rho||)$ and let $u$ be 
$(L,\rho)$-minimal; theorem \ref{periodic} yields the first equality below, formula (\ref{mean-action}) the second one:
$$\beta_{L-f}(\rho)=\frac{1}{|D|}
\int_D (L-f)(x,u,\nabla u)dx=
MA(L-f,T_u) . $$
Since $T_u\in H_r$, we get that
$$\beta_{L-f}(\rho)\ge \min_{T\in H_r}MA(L-f,T) . $$
To show the opposite inequality, we recall that 
$MA(L-f,\cdot)$ is affine, continuous and $H_r$ is the closed, convex hull of the currents $T_u$, with $u\in J_r (\rho)$; thus, it suffices
to prove that $MA(L-f,T_u)\ge\beta_{L-f}(\rho)$ for any 
$u\in J_r (\rho)$. But this follows immediately from theorem \ref{periodic}.

\qed

\begin{lemma}\label{continuity} Let $L$ satisfy hypotheses H1)-H4) of the introduction. Then the function 
$ T\longmapsto MA(L,T)$ is affine and continuous on $H_r$.

\end{lemma}

\proof We refer the reader to \cite{Bernard-Buffoni} for the proof that $MA(L,\cdot)$ is affine; we prove that it is continuous on 
$H_r$. 

Let $T\in H_r$ and let the measure 
$\mu_T$ be defined as in (\ref{mutti}); then  
\cite{Bernard-Buffoni} implies that there is a multi-vector field 
$X\in L^1(\mu_T)$ such that $T=X\wedge\mu_T$ and 
$X=(X_1,X_2,\dots,X_n,1)$ in the coordinates introduced above; moreover,
$$MA(L,T)=\int_{\T^n\times\T}L(x,u,X_1(x,u),\dots,X_n(x,u))d\mu_T(x,u) . $$
Let $\gamma_T$ be the push-forward of the measure $\mu_T$ by the map
$$(x,u)\longmapsto(x,u,X_1(x,u),\dots,X_n(x,u))  .  $$
By the formula above, we have that
$$MA(L,T)=\int_{\T^n\times\T\times\R^n}
L(x,u,p)d\gamma_T(x,u,p).$$
It is easy to see the following: if $T_k$ is a sequence in $H_r$, then it converges weak${}^\ast$ to $T$ if and only if the measures
$\gamma_{T_k}$ converge weak${}^\ast$ to $\gamma_T$. We also note that, by definition, the support of $\gamma_T$ with 
$T\in H_r$ is contained in $\T^n\times\T\times B(0,r)$; since on this set $L$ is bounded, we get that the linear function
$$\gamma_T\longmapsto
\int_{\T^n\times\T\times\R^n}
L(x,u,p)d\gamma_T(x,u,p)$$
is continuous; by the aforesaid this implies that also $MA(L,\cdot)$ is continuous.

\qed

We have shown that, generically, there is only one periodic minimizer, i. e. the first part of proposition \ref{generic_rho}; by the next lemma, this implies that all the
$(L,\rho)$-minimizers induce the same current, i. e. the second part of proposition \ref{generic_rho}.

\begin{lemma}\label{unicurrent} Let $\rho\in\Q^n$, and let us suppose that there is a unique periodic $(L,\rho)$-minimizer $u$; let us call $T_u$ the current it induces. Let $v$ be any $(L,\rho)$-minimizer; then 
$T_u=T_v$.

\end{lemma}

\proof By lemma 3.1 and integer translation, we can suppose that 
$u(0)\le v(0)<u(0)+1$. By \cite{Bangert}, there is a vector 
$\gamma\in\R^n$ such that
\begin{equation}
\lim_{t\to-\infty}
||v-u||_{C^1(\{ x\quad\colon\quad x\cdot\gamma<t \})}=0$$
$$\hbox{and}\quad
\lim_{t\to+\infty}||v-u-1||_{C^1(\{ x\quad\colon\quad x\cdot\gamma>t \})}=0 . 
\end{equation}
Let $\omega\in\Omega^{(0)}_n$; by the formula above, we can fix $t>0$ so large that
$$|\omega(x,v)\cdot(\nabla v,1)-
\omega(x,u)\cdot(\nabla u,1)|<\epsilon
\quad\hbox{if}\quad |x\cdot\gamma |>t  .  $$
For this $t$,
$$|T_v(\omega)-T_u(\omega)|\le
\lim_{R\to\infty}\frac{1}{|B(0,R)|}\int_{B(0,R)}
|\omega(x,v)\cdot(\nabla v,1)-\omega(x,u)\cdot(\nabla u,1)|dx=$$
$$\lim_{R\to\infty}\frac{1}{|B(0,R)|}
\int_{B(0,R)\cap\{ x\quad\colon\quad |x\cdot\gamma|>t \}}
|\omega(x,v)\cdot(\nabla v,1)-\omega(x,u)\cdot(\nabla u,1)|dx\le
\epsilon  .  $$
Since $\epsilon$ is arbitrary, the last formula implies the thesis.

\qed

\section{Uniqueness of the minimizing current within a given homology class: irrational case}\label{irrational currents}

\begin{lemma}\label{uniqueness}
Let \begin{itemize}
\item $L$ be a Lagrangian on $\R^{2n+1}$ satisfying Hypothesis (H1-4).
\item $\rho$ be a vector in  $\rho\in\R^n\setminus\Q^n$
\item $u_1,u_2$ be minimizers in  $M_\rho$.
\end{itemize}
Then  $T_{u_1}=T_{u_2}$.
\end{lemma}

\proof 
For $l=1,2$, we define
$$u^{\alpha^-}_l(x)=\sup\{ u_l(x+k)+j\quad\colon\quad \rho\cdot k+j<
\alpha \}  $$
and
$$u^{\alpha^+}_l(x)=\inf\{ u_l(x+k)+j\quad\colon\quad \rho\cdot k+j>
\alpha \} . $$
We recall a few results of Bangert's on the properties of 
$u_l^{\alpha^\pm}$. We set
$$\Gamma=\{ (k,j)\in\Z^n\times\Z\quad\colon\quad k\cdot\rho+j=0 \}
. $$

\noindent 1) It is proven in proposition 5.6 of \cite{Bangert} that 
$u_l^{\alpha^+}=u_l^{\alpha^-}$ save for at most countably many 
$\alpha$, for which $u_l^{\alpha^+}>u_l^{\alpha^-}$. We call 
$u_l^\alpha$ their common value, defined for $\alpha$ outside a countable set.

\noindent 2) By the same proposition, $u_l^{\alpha^\pm}$ is 
$\Gamma$-periodic; it follows from the definition that 
$u_l^{\alpha^-}\le u_l^{\alpha^+}$. Let $M$ be the projection of 
${\rm rat}(\rho,1)\subset\R^n\times\R$ on $\R^n$; corollary 4.6 of 
\cite{Bangert} implies that, for all $\epsilon>0$ we can find $C>0$ such that
$$u_l^{\alpha^+}(x+z)-u_l^{\alpha^-}(x+z)\le\epsilon
\quad\hbox{if}\quad x\in M,z\in M^\perp,\quad ||z||\ge C . $$

\noindent 3) By \cite{Bangert unique}, there is $a\in\R$ such that 
$u_1^\alpha=u_2^{\alpha+a}$. 

\noindent 4) Setting $\alpha_1=0$, $\alpha_2=-a$, we have by the last point and the definition of $u_l^\alpha$ that, for $l=1,2$,
$$u_1^{\alpha_l^-}\le u_l\le u_1^{\alpha_l^+} . $$
In the formula above, there are only two possibilities: either there are two equality signs, or there are two strict inequalities.

\vskip 2pc

We define $\tilde G$ as the closure of
$$\cup_{\alpha\in\R}
\{
(x,u^{\alpha^\pm}_1(x))\quad\colon\quad x\in\R^n
\}    $$
and we call $G$ the projection of $\tilde G$ on $\T^n\times\T$. 

\noindent{\bf Observation 1.} We assert that for $l=1,2$, 
$supp(T_{u_l})\subset G$. Indeed, let $\omega$ be a continuous 
$n$-form compactly supported on $(\T^n\times\T)\setminus G$; we shall show that $T_{u_i}(\omega)=0$. For starters, $\omega$ induces a periodic form 
$\tilde\omega$ on $\R^n\times\R$. Since $\omega$ is compactly supported on $\T^n\times\T$, the distance between the support of 
$\omega$ and $G$ is positive, which implies that the distance between the support of $\tilde\omega$ and $\tilde G$ is positive. By point 2) above, this implies that, for $C$ large enough,
$$
supp(\tilde\omega)\cap
\{
(x,x_{n+1})\in\R^n\times\R\quad\colon\quad 
u_1^{\alpha_l^-}(x)<x_{n+1}<u_1^{\alpha_l^+}(x)
\}
\subset $$
$$ M\times(B(0,C)\cap M^\perp) . 
$$
In particular, if $s$ is the dimension of $M$, we get that
$$|\{ 
x\in B(0,R)\quad\colon\quad (x,u_l(x))\in supp(\tilde\omega)
 \}|     \le C_1R^{s}  .  $$
Since $\rho$ is irrational, $s< n$, and thus
$$|T_{u_l}(\omega)|=\lim_{R\to\infty}
{1\over{|B(0,R)|}}
\left\vert
\int_{B(0,R)}\tilde\omega(x,u_l(x))\cdot(\nabla u_l,1)(x)dx
\right\vert   \le$$
$$\lim_{R\to\infty}
{{C_2R^{s}}\over{|B(0,R)|}} =0  $$
because $s<n$, since $\rho$ is irrational. 

\noindent{\bf Observation 2.} We assert that Lemma \ref{uniqueness}    follows if we prove that $\mu_{T_{u_1}}=\mu_{T_{u_2}}$. Indeed, let
$$X_l\colon \tilde G\cup
\{
(x,u_l(x))\quad\colon\quad x\in\R^n
\}
\to\Lambda_n(\R^{n+1})$$ 
be defined by
$$X_l(x,u_1^{\alpha\pm}(x))=(\nabla u_1^{\alpha\pm}(x),1),\qquad
X_l(x,u_{l}(x))=(\nabla u_{l}(x),1). $$
In the formula above, we have written the coordinates of $X_l$ with respect to  the basis $\{ e_i \}$ of $\Lambda_n(\R^{n+1})$ which is dual to the basis 
${\rm d}\hat x_i$ of 
$\Lambda^n(\R^{n+1})$ we defined in section \ref{current}.

By point 3) above, we have that $X_1=X_2$ on $G$; we call $X$ their common value on this set. Now $T_{u_l}$ is supported on $G$, where $X$ is defined; 
clearly, the observation follows if we prove that 
$T_{u_l}=X\wedge\mu_{T_{u_l}}$. 

To show this, we recall that, by 
(\ref{misura}),  
$\mu_{T_{u_l}}$ is the weak${}^\ast$ limit of the measures $\mu_{l,R}$ on $\T^n\times\T$ defined by
$$\int_{\T^n\times\T}
f(x,x_{n+1})d\mu_{l,R}(x,x_{n+1})=
{1\over{|B(0,R)|}}\int_{B(0,R)}f(x, u_l(x))dx$$
for all continuous functions $f$. Now, it follows from theorem 4.5 of \cite{Moser} that 
$X_l$ is Lipschitz on the union of  $G$ and  the graph of $u_l$; if 
$\tilde X_l$ is a Lipschitz extension of $X_l$ to
$\T^n\times\T$, we get  that
\begin{eqnarray*}
T_{u_l}(\omega) &=&
\lim_{R\to\infty}
{1\over{|B(0,R)|}}
\int_{B(0,R)}\omega(x,u_l(x))\cdot(\nabla u_l(x),1)dx \\
&=&\lim_{R\to\infty}
\int_{\T^n\times\T}\omega(x,x_{n+1})\cdot\tilde X_l(x,x_{n+1})
d\mu_{l,R}(x,x_{n+1})  \\
&=&\int_{\T^n\times\T}\omega(x,x_{n+1})\cdot\tilde X_l(x,x_{n+1})
d\mu_{T_{u_l}}(x,x_{n+1})  \\
&=&\int_{\T^n\times\T}\omega(x,x_{n+1})\cdot X(x,x_{n+1})
d\mu_{T_{u_l}}(x,x_{n+1})    
\end{eqnarray*}
where the first equality is the definition of $T_{u_l}$; the second one follows from the definition of $\mu_{l,R}$ and the fact that 
$\tilde X_l=(\nabla u_l,1)$ on the graph of $u_l$. The third equality follows since 
$\omega\cdot\tilde X_l$ is a continuous function on $\T^n\times\T$ and $\mu_{l,R}$ converges weakly. We note that, by (\ref{mutti}), if 
$T_{u_l}$ is supported on $G$, the measure $\mu_{T_{u_l}}$ is supported on $G$ too; since on this set $\tilde X_1=\tilde X_2=X$, the last equality follows.

The formula above implies that $T_{u_l}=X\wedge\mu_{T_{u_l}}$.

\noindent{\bf Observation 3.} We define the map
$$\tilde\Phi\colon\tilde G\to\R^{n+1},\qquad 
\tilde\Phi(x, u_1^{\alpha^\pm}(x))=(x,\rho\cdot x+\alpha) .  $$
We recall from \cite{Moser} that this map quotients to a map
$\Phi\colon{G}\to\T^{n+1}$.
We call $P$  the canonical projection 
$\T^{n}\times \T \longrightarrow \T^{n}$, i. e. 
$P(x,x_{n+1})=x$. We shall prove the following three facts.
\begin{itemize}
  
  \item  $P_\sharp(\mu_{T_{u_l}})$ and
  $(P\circ\Phi)_\sharp(\mu_{T_{u_l}})$ are the Lebesgue measure
  
  \item the measures $\mu_{T_{u_l}}$ on $G$  are invariant by the map
$$\psi_{k}\colon G\to G,\qquad \psi_k\colon (x,u^\alpha_1(x))\to
(x,u_1^{\alpha}(x+k))  $$
  
  \item the measures $\Phi_\sharp(\mu_{T_{u_l}})$ on $\T^n\times\T$ are invariant by the map 
$$\tilde\psi_k\colon\T^n\times\T\to\T^n\times\T,\quad
\tilde\psi_{k}(x,x_{n+1})=(x,x_{n+1}+k\cdot\rho) . $$
 
\end{itemize}

For the first statement, we note that, if $f\colon\T^n\to\R$ is continuous, then
\begin{eqnarray*}
\int_{\T^n}fdP_\sharp(\mu_{T_{u_l}}) &=&
\int_{\T^n\times\T}f(x)d\mu_{T_{u_l}}(x,x_{n+1}) \\
=\lim_{R\to\infty}{1\over{|B(0,R)|}}
\int_{B(0,R)}f(x)dx &=& \int_{\T^n}f(x)dx
\end{eqnarray*}
where the second equality comes from (\ref{misura}) and the last one from the periodicity of $f$. This proves that 
$P_\sharp(\mu_{T_{u_l}})$ is Lebesgue. The statement for 
$(P\circ\Phi)_\sharp(\mu_{T_{u_l}})$ follows as above, noting that 
$P\circ\Phi(x,u^{\alpha^\pm}(x))=x$.

The second statement follows from Lemma \ref{translation, current} and (\ref{mutti}).

To prove the third statement, we note that
$u_1^{\alpha+k\cdot\rho}=u_1^\alpha(x+k)$; we can rewrite this fact as
$$\Phi\circ\psi_k(x,u^{\alpha^\pm}(x))=
\tilde\psi_k\circ\Phi(x,u^{\alpha^\pm}(x)) . $$
This implies the first equality below, while the second one follows from the previous point.
$$(\tilde\psi_k)_\sharp\Phi_\sharp(\mu_{T_{u_l}})=
\Phi_\sharp(\psi_k)_\sharp(\mu_{T_{u_l}})=
\Phi_\sharp(\mu_{T_{u_l}})  .  $$
The last formula proves the invariance of 
$\Phi_\sharp(\mu_{T_{u_l}})$. 

\noindent{\bf Observation 4.} We assert that 
$\Phi_\sharp(\mu_{T_{u_l}})$ is the Lebesgue measure. We prove this using observation 3 and the Fourier transform; we set
$$m_{k,j}=\int_{\T^n\times\T}
e^{-2\pi i[k\cdot x+jx_{n+1}]}d\Phi_\sharp(\mu_{T_{u_l}})(x,x_{n+1}) . 
$$
We choose $\tilde k$ such that $\rho\cdot\tilde k\not\in\Q$; invariance under 
$\tilde\psi_{\tilde k}$ implies the second equality below
\begin{eqnarray*}
m_{k,j}&=&\int_{\T^n\times\T}e^{-2\pi i[k\cdot x+jx_{n+1}]}
d\Phi_\sharp(\mu_{T_{u_l}})(x,x_{n+1}) \\
&=&\int_{\T^n\times\T}e^{-2\pi i[k\cdot x+jx_{n+1}]}
d(\tilde\psi_{\tilde k})_\sharp\Phi_\sharp(\mu_{T_{u_l}})(x,x_{n+1})
\\
&=&\int_{\T^n\times\T}e^{-2\pi i[k\cdot x+j(x_{n+1}+\tilde k\cdot\rho)]}
d\Phi_\sharp(\mu_{T_{u_l}})(x,x_{n+1})=
e^{-2\pi ij(\tilde k\cdot \rho)}m_{k,j} . 
\end{eqnarray*}
Since $\rho\cdot\tilde k$ is irrational, we deduce that $m_{k,j}=0$ unless $j=0$. Since the marginal of $\Phi_\sharp(\mu_{T_{u_l}})$ on $\T^n$ is the Lebesgue measure by the first point of observation 3, we see that 
$m_{k,0}=0$ unless $k=0$; in other words, 
$\Phi_\sharp(\mu_{T_{u_l}})$ has the same Fourier transform as the Lebesgue measure, which implies that it is the Lebesgue measure.

\noindent{\bf End of the proof.} We prove that 
$\mu_{T_{u_1}}=\mu_{T_{u_2}}$; by observation 2, this implies the thesis.

If $\Phi$ were injective, observation 4 would imply that 
$\mu_{T_{u_1}}$ and $\mu_{T_{u_2}}$ coincide. Indeed, 
\begin{eqnarray*}
\mu_{T_{u_1}}(A)&=&\mu_{T_{u_1}}(\Phi^{-1}(\Phi(A)))=
\Phi_\sharp\mu_{T_{\mu_1}}(\Phi(A))={\mathcal L}^{n+1}(\Phi(A)) \\
&=&\Phi_\sharp\mu_{T_{\mu_2}}(\Phi(A))=
\mu_{T_{u_2}}(\Phi^{-1}(\Phi(A)))=\mu_{T_{u_2}}(A)
\end{eqnarray*}
where the third and fourth equalities come from observation 4.
The same argument would apply if we could prove that the set on which $\Phi$ is not injective is negligible for 
$\mu_{T_{u_1}}$ and $\mu_{T_{u_2}}$.  But we saw in point 4) at the beginning of the proof that the set on which $\Phi$ is two to one is exactly the union of the boundaries of the  gaps of $G$, 
which have the form 
$$\{ (x,u_1^{\alpha^\pm}(x))\quad\colon\quad x\in\R^n \}$$
projected to the torus. By countable additivity, it suffices to prove that each piece
$$\{ (x,u_1^{\alpha^\pm}(x))\quad\colon\quad x\in [0,1]^n \}$$
has measure zero. But the measures $\mu_{T_{u_i}}$ are invariant by the action of $\psi_{k}$; thus, if one of the sets above had positive measure, the measure of the whole torus would be infinite, a contradiction.

 \qed

\section{Proof of the main theorem}\label{last}

\begin{theorem}\label{main}
Let $L$ be a Lagrangian on $\R^{2n+1}$ satisfying Hypothesis (H1-4). Then there exists a residual subset $\mathcal{O}(L)$ of $C^{\infty}(\T^{n+1})$ such that for any 
$f \in \mathcal{O}(L)$, 
\begin{itemize}
  \item for any $\rho \in\R^n $, all the $(L-f,\rho)$-minimizers induce the same current; if $\rho \in\Q^n $, there is a unique periodic 
  $(L,\rho)$-minimizer
  \item for any $c \in\R^n$, all the $(L-f-c)$-minimizers induce the same current
  \item there exists an open dense subset $U(L,f)$ of $\R^n$, such that for any $c \in U(L,f)$, we have $\rho :=\alpha_{L-f}'(c) \in \Q^n$.
 
\end{itemize}
\end{theorem}

\proof

\noindent{\bf First statement.} 
Set
	\[ \mathcal{O}(L) := \bigcap_{\rho \in \Q^n} \mathcal{O}(L,\rho)
\]
where $\mathcal{O}(L,\rho)$ comes from Proposition \ref{generic_rho}. Then $\mathcal{O}(L)$ is residual in \newline
$C^{l,\gamma}(\T^{n+1})$. Take $f\in\mathcal{O}(L)$. 
We remark that, if $\rho \in \R^n$ is irrational, then by Lemma
\ref{uniqueness}, all the  $(L-f,\alpha_{L-f}' (c))$-minimizers induce the same current. On the other hand, if $\rho \in \Q^n$, then by the definition of $\mathcal{O}(L)$ there is only one periodic $(L-f,\rho)$- minimizer; moreover, by lemma \ref{unicurrent}, the $(L,\rho)$-minimizers induce a unique current. This proves the first part of the theorem. 

\vskip 1pc

\noindent{\bf Second statement.}  Take  $f \in \mathcal{O}(L)$ and 
$c \in H^n (\T^{n+1})$. Then, if $\alpha_{L-f}' (c) \not\in \Q^n$, we know by Lemma \ref{uniqueness} that there exists a unique 
$(L-f,\alpha_{L-f}' (c))$-minimizing current, hence  there exists a unique $(L-f-c)$-minimizing current. If $\alpha_{L-f}' (c) \in \Q^n$, by the definition of 
$\mathcal{O}(L)$,  all the  $(L-f,\alpha_{L-f}' (c))$-minimizers induce the same current, hence  there exists a unique $(L-f-c)$-minimizing current. This proves the second part of the theorem.

\vskip 1pc

\noindent{\bf Third statement.} By the first statement above and Theorem \ref{Senn2}, if $f \in \mathcal{O}(L)$ and $\rho \in \Q^n$, the dimension of $D_{\rho}(L-f)$ is $n$. Now $D_\rho(L-f)\subset\R^{n+1}$; let 
$\mbox{P}\colon\R^{n+1}\to\R^n$ be the projection to the first $n$ coordinates, and let $\mbox{int}(D_{\rho}(L-f))$ denote the interior of $D_{\rho}(L-f)$. Since the dimension of $D_{\rho}(L-f)$ is $n$, and this set is not vertical, we easily get that 
$\mbox{P}\{ [ \mbox{int}(D_{\rho}(L-f)) ] \}$ is an open set.

Set 
	\[ U(L,f) := \bigcup_{\rho \in \Q^n} 
	\mbox{P}\{ [ \mbox{int}(D_{\rho}(L-f)) ] \},
\]
then $U(L,f)$ is open  in $\R^n$, and it is dense in $\R^n$ by Proposition \ref{dense}. Besides, if $\rho \in \Q^n$ and 
$c \in \mbox{P}\{ [ \mbox{int}(D_{\rho}(L-f)) ] \}$, then by Proposition \ref{generic_rho} there is a unique  periodic $(L-f-c)$-minimizer with slope $\alpha_{L-f}'(c)$.
\qed
\appendix
\section{A bit of topology}
We denote by 
\begin{itemize}
	\item  $B(0,r)$ the closed ball in $\R^{n}$ of radius $r$, centered at the origin
	\item  $ \langle .,. \rangle $ the canonical inner product in $R^n$.
	\end{itemize}

\begin{lemma}\label{appendix_lm1}
Let $f$ be a continuous map from $\R^{n}$ to itself, such that 
  for any $x$ in $\R^n \setminus \{0\}$, we have $\left\langle x, f(x) \right\rangle > 0$.
 Then for any neighborhood $U$ of zero in $\R^{n}$, $f(U)$ contains a neighborhood of zero.
\end{lemma}

\proof

By modifying $f$ outside some neighborhood of zero, we may assume that  $\left\|f(x)\right\|$ goes to infinity when $\left\|x\right\|$ goes to infinity.
Thus, setting $\tilde{f}(\infty):=\infty$, $f$ extends to a continuous self-map $\tilde{f}$ of 
$\R^n\cup\{ \infty \}$, the one-point compactification of $\R^n$. We identify $\R^n\cup\{ \infty \}$ with $\Sbar^n$ by the stereographic projection, i. e. by the map
$$\psi\colon\{ (x,z)\quad\colon\quad x\in\R^n,z\in\R,\quad 
|x|^2+z^2=1 \}\to\R^n$$
defined by
$$\psi(x,z)={1\over{1-z}}x  . $$
We consider the continuous map $F\colon \Sbar^n\to \Sbar^n$ given by 
$F(x,z)=\psi^{-1}\circ\tilde f\circ\psi(x,z)$.  Next, we observe that
$\langle\psi(x,z),\psi(-x,-z)\rangle<0$ save when $x=0$; in other words, if two points on $S^n$ are diametrically opposite, then the internal product of their $\psi$-images is negative. In particular, if 
$(x,z)$ and $F(x,z)$ were diametrically opposite, then we would have that $\langle\psi(x,z),\tilde f(\psi(x,z))\rangle<0$; but this is excluded by our hypotheses on $f$. Since 
$(x,z)$ and $F(x,z)$ are never diametrically opposite,  $F$  is homotopic to the identity by the shortest geodesic homotopy. Therefore  $F$ has degree one, hence it is onto; as a consequence, $\tilde f$ is onto too.

Now we want to show that for any $\delta >0$, there exists $\epsilon >0$ such that $B(0,\epsilon) \subset f(B(0,\delta))$. Assume not. Then for any positive integer $k$, there exists $x_k$, such that $\left\|x_k \right\| \leq\frac{1}{k}$ and $x_k$ is not in $f(B(0,\delta))$. Since $f$ is onto, there exists $y_k$, such that $\left\|y_k \right\| \geq \delta$ and $f(y_k)=x_k$. Take a limit point $y$ of $y_k$ in the $n$-dimensional sphere. Since $||y_k||\ge\delta$, we have $y\not=0$. But, since $f$ is continuous, we get $f(y)=0$,
a contradiction with  $\left\langle y, f(y) \right\rangle > 0$.
\qed

\section{A bit of linear algebra}
We say an affine subspace of $\R^n$ is rational if it is defined by equations of the form $\left\langle c_i, h \right\rangle =  \tau_i$, $i=1, \ldots s$, where  $c_i$, $i=1, \ldots s$, are integer vectors, and 
$\tau_i \in \Z$, $i=1, \ldots s$. The intersection of two rational affine subspaces is a rational affine subspace, so given $\rho \in \R^n$, there exists a smallest rational affine subspace containing $\rho$. We denote it $A(\rho)$. 

With an abuse of notation, we shall set 
$$A(\rho)^\perp=\mbox{Vect}(c_1,\ldots c_s) . $$
In other words, $A(\rho)^\perp$ is the vector space orthogonal to 
$L(A(\rho))$, i. e. to the smallest space containing the differences 
$a-b$ with $a,b\in A(\rho)$.

We define $\rat(\rho,1)$ as the subspace of $\R^n$ generated by 
$\Z^n\times\Z\cap(\rho,1)^\perp$; we also define $M(\rho)$ as the projection on $\R^n$ of $\mbox{rat}(\rho,1)$.
Recall from \cite{vmb2} (Proposition 18) that the irrationality $I_{\Z}(\rho)$ of $\rho$ is the dimension of $A(\rho)$. The next lemma implies that
$$
I_{\Z}(\rho) = n-\dim M(\rho) = n -\dim \mbox{rat}(\rho,1).
$$

\begin{lemma}
 For any $\rho \in \R^n$, $M(\rho)$ is the vector subspace $A(\rho)^{\perp}$ of $\R^n$ orthogonal to  $A(\rho)$.
\end{lemma}
\proof
Assume $A(\rho)$ is defined by the equations $c_i\cdot v= \tau_i$, with $c_i \in \Z^n$, $\tau_i \in \Z$, $i=1,\ldots k$. Then $A(\rho)^{\perp}=\mbox{Vect}(c_1,\ldots c_k)$. 
Recall from \cite{Bessi09} that $M(\rho)$ is generated by the vectors $k \in \Z^n$ such that $\rho\cdot k \in \Z$. Thus $c_i \in M(\rho)$, $i=1,\ldots k$, whence $A(\rho)^{\perp} \subset M(\rho)$.

On the other hand, if $k \in \Z^n$  is such that $\rho\cdot k \in \Z$, the equations $c_i\cdot v= \tau_i$, $i=1,\ldots s$ together with
 $k\cdot v =\rho\cdot k $ define  a rational affine subspace $B$ of $\R^n$, containing $\rho$, and contained in $A(\rho)$, so by the definition of $A(\rho)$, $B=A(\rho)$. Therefore $k \in  \mbox{Vect}(c_1,\ldots c_s)$. Since $A(\rho)^{\perp}=\mbox{Vect}(c_1,\ldots c_s)$, we conclude that $A(\rho)^{\perp} \supset M(\rho)$.

\qed

\begin{lemma}\label{B2}
For any $\rho \in \R^n$, $\rho=\rho_1 + \rho_2$, where $\rho_1 \in  M(\rho) \cap \Q^n$, $\rho_2 \in  M(\rho)^{\perp}$, and $\rho_2$ is completely irrational  in  $M(\rho)^{\perp}$, that is, $\rho_2$ is not contained in any proper rational affine subspace of $M(\rho)^{\perp}$.
\end{lemma}
\proof

Let $P$ denote the orthogonal projection on $M(\rho)$; we begin to prove that $P(\rho)$ is rational. Let $w_1, \ldots w_k$ be vectors in $\Z^n$ which form a basis of $M(\rho) $; by the definition of $M(\rho)$ we get that $\rho\cdot w_i\in\Z$. Since
$P(\rho)\cdot w_i=\rho\cdot w_i$, we get that 
$P(\rho)\cdot w_i\in\Z$; if we set 
$$P(\rho)=a_1w_1+\dots a_kw_k, \quad
a=(a_1,\dots,a_k), \quad
b=(P(\rho)\cdot w_1,\dots,P(\rho)\cdot w_k)$$ 
and we define $W$ to be the matrix of the internal products $w_i\cdot w_j$, we see that
$Wa=b$, i. e. $a=W^{-1}b$. Since $W$ and $b$ have integer entries, this implies that $a$ is rational, which in turn implies that $P(\rho)$ is rational.

Now set $\rho_2 = \rho -P(\rho)$, we have $\rho_2 \in M(\rho)^{\perp}$. Assume $\rho_2$ lies in a rational affine subspace $B$ contained in $M(\rho)^{\perp}$. Then $P(\rho) +B$ is a rational affine subspace of $\R^n$, and it contains $\rho$, so it contains $A(\rho)$. On the other hand, the dimension of $P(\rho) +B$ is at most 
$\dim  M(\rho)^{\perp} = \dim A(\rho)$, so $B=M(\rho)^{\perp}$. Thus $\rho_2$ is completely irrational  in  $M(\rho)^{\perp}$.
\qed

  {\small

\bigskip

\noindent
Dipartimento di Matematica, Universit\`a Roma Tre, Italy\\
e-mail : bessi@mat.uniroma3.it\\
D\'epartement de Math\'ematiques, Universit\'e Montpellier 2, France\\
e-mail : massart@math.univ-montp2.fr
}

\end{document}